\documentclass[12pt]{article}

\usepackage{multirow}

\usepackage{ifthen}
\newboolean{PrintVersion}
\setboolean{PrintVersion}{false} 
\usepackage{amsmath,amssymb}
\usepackage{amsthm}
\usepackage{mathtools}
\usepackage{graphicx}
\usepackage{caption}
\usepackage{bbm}
\usepackage{color}
\usepackage{dcolumn}
\usepackage{bm}
\usepackage{float}
\usepackage{hyperref} 
\usepackage{apacite}
\definecolor{background-color}{gray}{0.98}
\definecolor{steelblue}{rgb}{0.27, 0.51, 0.71}
\definecolor{brickred}{rgb}{0.8, 0.25, 0.33}
\definecolor{bluegray}{rgb}{0.4, 0.6, 0.8}
\definecolor{amethyst}{rgb}{0.6, 0.4, 0.8}

\hypersetup{
    plainpages=false,       
    pdfpagelabels=true,     
    bookmarks=true,         
    unicode=false,          
    pdftoolbar=true,        
    pdfmenubar=true,        
    pdffitwindow=false,     
    pdfstartview={FitH},    
    pdfnewwindow=true,      
    colorlinks=true,        
    linkcolor=steelblue,         
    citecolor=bluegray,        
    filecolor=magenta,      
    urlcolor=cyan           
}
\ifthenelse{\boolean{PrintVersion}}{   
\hypersetup{	
    citecolor=black,%
    filecolor=black,%
    linkcolor=black,%
    urlcolor=black}
}{} 

\theoremstyle{plain}
\newtheorem{theorem}{Theorem}[section]
\newtheorem{lemma}[theorem]{Lemma}

\newtheorem{proposition}[theorem]{Proposition}

\theoremstyle{definition}

\newenvironment{keywords}
  {\noindent \textbf{Keywords: }}

\newcommand{\Nset}[1]{\mathcal{N}_{#1}}
\newcommand{\set}[1]{\mathcal{#1}}
\newcommand{\event}[1]{\mathcal{#1}}

\newcommand{\field}[1]{\mathbb{#1}}
\newcommand{\Reals}{\field{R}}

\newcommand{\Naturals}{\field{N}}
\newcommand{\NaturalsZero}{\field{N}_0}

\newcommand{\follows}{\sim}

\newcommand{\ffact}[2]{[{#1}]_{#2}}

\title{Bernoulli Sums: The only random variables that count}

\author{Pavel Shuldiner\thanks{We acknowledge the support of the Natural Sciences and Engineering Research Council of Canada (NSERC), PGSD2 - 535019 - 2019.} ~and R.W. Oldford\\
  University of Waterloo }

\date{\today}

\begin{document}
\maketitle
\begin{abstract}
A novel  multinomial theorem for commutative idempotents is shown to lead to new results about the moments, central moments, factorial moments, and their generating functions for any random variable $X = \sum_{i} Y_i $ expressible as a sum of Bernoulli indicator random variables $Y_i$.
The resulting expressions are functions of the expectation of products of the Bernoulli indicator random variables.  These results are used to derive novel expressions for the various moments in a large number of familiar examples and classic problems including: the binomial, hypergeometric, Poisson limit of the binomial, Poisson, Conway-Maxwell-Poisson binomial,  Ideal Soliton, and Benford distributions as well as the empty urns problem and the matching problem.  Other general results include expressions for the moments of an arbitrary count random variable in terms of its upper tail probabilities.
\end{abstract}
\begin{keywords}
factorial moments, central moments, binomial distribution, hypergeometric,  Poisson, the empty urns problem, the matching problem, CMP-binomial, Benford distribution, ideal Soliton, Stirling numbers, Bell numbers
\end{keywords}

\section{Introduction}
\label{sec:intro}
Consider the random variable
\[ X =  \sum_{i \in \set{I}} Y_i \]
which sums (not necessarily independent) Bernoulli random variables $Y_i  \in \{0, 1\}$ over some countable index set $\set{I}$.
When $Y_i$ is an indicator function for some event $\event{A}$, the \textit{Bernoulli sum} $X$ counts the number of occurrences of the event in the set $\set{I}$ and, as such,  arises in numerous applications of probability.
Of interest here is the determination of the moments, central moments, and factorial moments of any arbitrary Bernoulli sum.  

We develop expressions for these moments in terms of the expectation of products of the Bernoulli $Y_i$s.  This leads to novel proofs and/or  expressions for the moments in many well known problems and to novel approaches to determining such moments for any random variable expressible as a Bernoulli sum.

The paper is organized as follows.  
Section \ref{sec:preliminaries} shows that the power of a sum of idempotents is expressible in terms of the number of surjections from one finite set to another times a sum of their products.
This follows as a special case of the multinomial theorem.  Section \ref{sec:moments} builds on this to develop the main general results for the moments of a Bernoulli sum.  Both finite and infinite sums are considered and special attention is given to factorial moments and generating functions.  

These results are then applied to develop expressions for various classic distributions and problems in Section \ref{sec:classics}.  These include the binomial, poisson binomial, hypergeometric,  and Conway-Maxwell-Poisson binomial distributions, the Poisson limit of a binomial by moment convergence, and the classic empty  urns problem and the matching problem.   

Section \ref{sec:counts_generally} considers the moments for any count random variable, developing expressions based on the upper tail probability of that count.  
This general theory is then demonstrated on the geometric,  Poisson, Ideal Soliton, and Benford distributions.  The paper ends with a few concluding remarks in its last section.\section{Preliminaries}
\label{sec:preliminaries}
First, let $\Nset{k}$ denote the set $\left\{1, \ldots, k \right\}$ for any finite integer $k$, and $S(k, m)$ denote the number of surjections from $\Nset{k}$ to $\Nset{m}$.   If $k < m$, no surjective function exists and $S(k,m) = 0$; otherwise,  $S(k,m)$ can be written as
 \[S(k,m) = \sum_{v = 0}^{m-1}(-1)^{v} \binom{m}{v} (m-v)^{k}\]
 \cite<e.g., see> {wilf2005generatingfunctionology}.
 The number $S(k, m)$ figures prominently in the closed form expressions which follow.  
 
 In particular, $S(k,m)$ (for all $m \le k$) will be shown to appear in expressions for the $k$th moments of a Bernoulli sum.  To calculate the smaller moments of importance in statistical inference (say $k \le 4$), it will be convenient therefore to have $S(k,m)$ evaluated for a few $m \le k$.   Whenever $k$ is at least as large as the second argument, the following values obtain:
 $S(k,0) = 0$, $S(k, 1) = 1$, $S(k,2) = 2^k -2$, $S(k, 3) = 3^k - 3 \cdot 2^k +3$, and $S(k, 4) = 4^k - 4 \cdot 3^k + 6 \cdot 2^k -4$.  These will appear in calculations up to the $4$th moment (e.g. to determine kurtosis). Again, note that 
 $S(k,m) = 0$ whenever $m > k$.
 
Moments of $X$ are expectations of powers of $X$ which, in the case of $X = \sum_{i = 1}^n Y_i$,   suggests beginning with a multinomial theorem, here expressed fairly generally as:

\begin{theorem}
\label{th:classic_multinomial_theorem}
\textbf{The multinomial theorem.} 
Let $y_1,\ldots, y_n$ be a sequence of commutative elements over some ring and fix $k \ge 1$. Then 
$$(y_1 + \dotsb + y_n)^k  = \sum_{\substack{\ell_1 + \dotsb + \ell_n = k\\ \ell_i \ge 0, \forall i}}\binom{k}{\ell_1,\dotsb, \ell_n} ~ y_1^{\ell_1}\dotsb y_n^{\ell_n}.$$
\end{theorem}

\begin{proof}
E.g., see \citeA{goulden2004combinatorial} or \citeA{wilf2005generatingfunctionology}.
\end{proof}


When the $y_i$s are also idempotents, as they are in the definition of a Bernoulli sum,  
only those $y_i$ with $\ell_i \ge 1$  remain and simplify to $y_i^{\ell_i} = y_i$.
This leads to the following version of a 
multinomial theorem where now $S(k,m)$ appears.

\begin{proposition}
\label{prop:multinom_specialization}
	Let $y_1,\ldots, y_n$ be a sequence of commutative idempotents over some ring. Then 

	$$(y_1 + \dotsb + y_n)^{k} = 
	\sum_{m=1}^{k} S(k,m) \sum_{\substack{{\{i_1, \ldots,  i_m \}\subseteq \Nset{n}}}} y_{i_1}\dotsb ~y_{i_m}$$
where $S(k,m)$ is the number of surjections from $\Nset{k}$ onto $\Nset{m}$. \\
(It is understood that the interior sum has $m \leq \min{\left\{k, n\right\}}$.) 
\end{proposition}

\begin{proof} 
Naively
expanding 
${(y_1  + \dotsb + y_n)}^k$ gives
\begin{equation}
\label{noncommutative_expansion}
(y_1 + y_2 + \dotsb + y_n)^k = \sum_{(j_1, j_2,\ldots, j_k)\in \Nset{n}^k} y_{j_1}y_{j_2}\dotsb y_{j_k}.	
\end{equation}
Let $\mathcal{F}$ denote the set of all functions $f:\Nset{k}\to\Nset{n}$. For a product $y_{j_1}y_{j_2}\dotsb y_{j_k}$
 on the right-hand side of Equation \ref{noncommutative_expansion}, let $f$ be the function which maps $\ell \in \Nset{k}$ to $j_\ell$. Since every $j_\ell \in \Nset{n}$, this defines a function $f \in \mathcal{F}$. 
 Conversely, a unique summand of the form $y_{f(1)}y_{f(2)}\dotsb y_{f(k)}$ can be assigned to each function $f \in \mathcal{F}$.
 That is, the naive expansion of $(y_1 + \dotsb + y_n)^k$ 
 results
in $n^k$ summands of the form $y_{f(1)}y_{f(2)}\dotsb y_{f(k)}$ for some function $f:\Nset{k} \to \Nset{n}$.
 
Since
$y_1,\ldots, y_n$ are commutative idempotents,  each product  $y_{j_1} \dotsb ~y_{j_k}$ resolves to a unique  $y_{i_1} \dotsb ~y_{i_m}$ with indices $i_1 < \ldots < i_m$, for some $m \in  \Nset{k}$. 
Equation \ref{noncommutative_expansion} then becomes 
\begin{equation}
(y_1 + y_2 + \dotsb + y_n)^k = \sum_{(j_1, j_2,\ldots, j_k)\in \Nset{n}^k} y_{j_1}y_{j_2}\dotsb y_{j_k} = \sum_{m=1}^{k}    a(k,m)  
	     \sum_{\substack{
	                        {\{i_1,\ldots, i_m\}\subseteq \Nset{n}} \\
	                        }
	                        }
                y_{i_1}\dotsb ~y_{i_m}.
\end{equation}
Here $a(k,m)$ counts the number terms $y_{i_1}\dotsb~ y_{i_m}$ that simplify to $y_{j_1}\dotsb y_{j_k}$.
It remains only to show that $a(k,m)$ equals $S(k, m)$,  the number of surjective maps from $\Nset{k}$ onto $\Nset{m}$.
 
To see this, first fix $\{i_1, \ldots,  i_m \}\subseteq \Nset{n}$ and let $F \subseteq \mathcal{F}$ denote the 
subset of 
functions 
for which $y_{f(1)}\dotsb ~ y_{f(k)}$ simplifies to $y_{i_1}\dotsb y_{i_m}$. 
The count $a(k,m)$ is identical to $|F|$.
Then consider the set, $G$, 
of all surjections ${g: \Nset{k} \to \{i_1,\ldots, i_m\}}$, which must have size $|G| = S(k,m)$.
If $F = G$, then $ |F| = |G|$ and $a(k,m) = S(k,m)$, as required.

 Now $F = G$ iff every $f \in F$ is also in $G$ and every $g \in G$ is also in $F$.
 If $f \in F$, then $y_{f(1)}\dotsb ~ y_{f(k)} = y_{i_1}\dotsb y_{i_m}$,
 giving $f(\Nset{k}) = \{i_1,\ldots, i_m\}$, and hence $f \in G$.
 If $g \in G$, then clearly 
 \[ y_{g(1)}\dotsb ~ y_{g(k)} = \prod_{\ell \in g(\Nset{k})} y_{\ell} = \prod_{\ell \in \{i_1,\ldots, i_m\}} y_{\ell},\]
and so $g\in F$. 
\end{proof}
Note that the inner sum $\sum_{\substack{{\{i_1, \ldots,  i_m \}\subseteq \Nset{n}}}} y_{i_1}\dotsb ~y_{i_m}$ vanishes whenever $m > n$ and hence $(y_1+\dotsb + y_n)^k$ is expressible as a sum of at most $min(k,n)$ terms involving the coefficients $S(k,m)$.


A generalization of the result to powers of infinite sums, subject to
convergence having been settled for any particular values of the $y_i$s,
is relatively straightforward.
\begin{proposition}
\label{prop:inf_power_expansion}
	Let $(y_i)_{i\ge 1}$ be a sequence of formal, commutative, idempotents over some ring. Then 

	$$\left(\sum_{i = 1}^{\infty}y_i\right)^{k} = 
	\sum_{m=1}^{k} S(k,m) 
	\sum_{\{i_1, \ldots, i_m\}\subset \Naturals } 
	y_{i_1}\dotsb y_{i_m}$$
	
	where $S(k,m)$ is the number of surjections from $\{1,\ldots, k\}$ onto $\{1,\ldots, m\}$. 
\end{proposition}

\begin{proof}
Since the $(y_i)_{i\ge 1}$ are formal, commutative idempotents over some ring, the proof 
follows that of Proposition \ref{prop:multinom_specialization}. 
\end{proof}



\section{Moments of Bernoulli sums}
\label{sec:moments}
Consider the Bernoulli sum random variable $X$ of Section \ref{sec:intro} with finite index set $\set{I}$ of size  $n$.
Set $\set{I}$ can always be re-indexed to have $X$ appear as
\[ X =  \sum_{i = 1}^{n} Y_i . \]
Expressions for the moments of $X$ can now be derived via Proposition \ref{prop:multinom_specialization}.
\begin{proposition}
\label{prop:bern_moments}
	When $X$ is expressible as a finite Bernoulli sum $X =  \sum_{i = 1}^{n} Y_i$, the $k$th \textbf{moment} of $X$ is expressible as 
\[E(X^k) = 
	\sum_{m=1}^{k} S(k,m) \sum_
	{\{i_1, \ldots,  i_m \}\subseteq \Nset{n}} 
	E(Y_{i_1}\dotsb Y_{i_m}).\]
\end{proposition}

\begin{proof}
Since $Y_i^2 = Y_i ~ \forall ~i$, it follows from Proposition \ref{prop:multinom_specialization} that
\[ X^k = (Y_1+\dotsb + Y_n)^k = 
	\sum_{m=1}^{k} S(k,m) 
	\sum_{{\{i_1, \ldots,  i_m \}\subseteq \Nset{n}}} 
	Y_{i_1}\dotsb Y_{i_m}.\]
The result follows by applying expectation $E(\cdot)$ (or, more generally, any linear) operator to each side.
\end{proof}
Since 
\[ E(Y_{i_1}\dotsb Y_{i_m}) = Pr(Y_{i_1} = 1, \dotsb, Y_{i_m} = 1)
\]
Proposition \ref{prop:bern_moments}  shows that the moments of \textit{any} finite Bernoulli sum random variable can be investigated via the joint distribution of those Bernoulli random variables used to construct it -- indeed, 
Proposition \ref{prop:bern_moments} could be rewritten in terms of this probability.

The  \textit{central moments} are generally of more statistical interest and a similar result is found for them by applying Proposition \ref{prop:bern_moments}.  In this case, let
$p_i = Pr(Y_i =1) = E(Y_i)$ denote the $i$th marginal mean in the sum and 
$\mu = E(X) = \sum_{i = 1}^n p_i$ the mean of $X$.
A similar expression for the $k$th central moment is given in
Proposition \ref{prop:bern_central_moments}.

\begin{proposition}
	\label{prop:bern_central_moments}
	When $X$ is expressible as a finite Bernoulli sum $X =  \sum_{i = 1}^{n} Y_i$, with $p_i = Pr(Y_i =1) = E(Y_i)$, then the $k$th \textbf{central moment} of $X$ is expressible as

\[
E((X-\mu)^k) = \left(-\mu\right)^{k} +
\sum_{\ell=1}^{k}
       \binom{k}{\ell} 
         \left(-\mu\right)^{k-\ell}
        \sum_{m=1}^{\ell} 
         S(\ell,m)
        \sum_{{\{i_1, \ldots,  i_m \}\subseteq \Nset{n}}}
         E(Y_{i_1}\dotsb Y_{i_m})
\]
where $\mu = E(X) = \sum_{i = 1}^n p_i$.
\end{proposition}

\begin{proof}
Applying the binomial expansion, then Proposition \ref{prop:bern_moments}, yields
\begin{align*}
E((X-\mu)^{k}) 
&= \sum_{\ell=0}^{k}
      \binom{k}{\ell} E(X^\ell)(-\mu)^{k-\ell}\\
&= \binom{k}{0}\cdot 1\cdot (-\mu)^{k} +
\sum_{\ell=1}^{k}
       \binom{k}{\ell} 
         (-\mu)^{k-\ell}
        E(X^{\ell})
         \\
&= (-\mu)^{k} +
\sum_{\ell=1}^{k}
       \binom{k}{\ell} 
         (-\mu)^{k-\ell}
        \sum_{m=1}^{\ell} 
         S(\ell,m)
        \sum_{{\{i_1, \ldots,  i_m \}\subseteq \Nset{n}}}
         E(Y_{i_1}\dotsb Y_{i_m})
         \\
\end{align*}
\end{proof}

Of course, whenever the $Y_i$s are also \textbf{independently distributed}, the above moment expressions (and those which follow) simplify by replacing $E(Y_{i_1}\dotsb Y_{i_m})$ by $p_{i_1} \dotsb p_{i_m}$, where each $Y_i \follows$ Bernoulli($p_i$).

\subsection{Moments of an infinite sequence}
\label{sec:moments_infinite}

Consider now an infinite sequence 
\[(Y_i)_{i \ge 1} = Y_1, Y_2, \ldots \]
of Bernoulli random variables and their sum
\[ X = \sum_{i = 1}^{\infty} Y_i \]
being such that $Pr(X < \infty) = 1$.  
(This condition is satisfied, for example, whenever the first moment of $X$ is bounded, that is, whenever $E(X) = \sum_{i=1}^{\infty}E(Y_i) = \sum_{i\ge 1} p_i  = \mu <\infty$. 
From this it follows (e.g., by the Borel-Cantelli lemma) that $Pr(\limsup_{n\to\infty} Y_n = 1) = 0$, and, so, that the
probability is zero that infinitely many of the $Y_i$s will be 1.)
%

In this case, Proposition \ref{prop:inf_power_expansion} gives the $k$th moment for this sum of countably infinite Bernoulli random variables (whenever all relevant sums converge).
\begin{proposition}
\label{prop:moment_inf_expansion}
Let $X = \sum_{i = 1}^{\infty} Y_i$ be the sum of the sequence of Bernoulli random variables $(Y_i)_{i \ge 1}$,
with $p_i = Pr(Y_i = 1) =  1 - Pr(Y_i =0)$,  then we may write
\[E(X^k) = 
	\sum_{m=1}^{k}
	 S(k,m)
	\sum_{{\{i_1, \ldots, i_m\}\subset \Naturals }} 
	~ E(Y_{i_1}\dotsb Y_{i_m}).\]
\end{proposition}

In the special case where the $Y_i$s are also independent, then Proposition \ref{prop:moment_inf_expansion} allows us to draw the interesting conclusion that \textit{a bounded first moment of $X$ implies that all higher order moments are also bounded}.  This result is formally given in 
Proposition \ref{prop:bounded_moments}:

\begin{proposition}
\label{prop:bounded_moments}

Let $(Y_i)_{i\ge 1}$ be a sequence of \textbf{independent} Bernoulli$(p_i)$ random variables 
with $\sum_{i\ge 1} p_i = \mu < \infty$. 
For the  Bernoulli (infinite) sum random variable $X = \sum_{i\ge 1} Y_i$, and for any $k \ge 1$, 
\[
E(X^k)  < \infty.
\]

\end{proposition}

\begin{proof}
By Proposition \ref{prop:moment_inf_expansion} and independence of the $Y_i$, 
\begin{align*}
E(X^k) 
  &= \sum_{m=1}^{k}
           S(k,m)
       \sum_{\{i_1, \ldots, i_m\}\subset \Naturals } 
           p_{i_1} \dotsb p_{i_m}\\
  &\le \sum_{m=1}^{k}
           S(k,m)
       \left( \sum_{i\ge 1} p_i \right)^m\\
&= \sum_{m=1}^{k}
    S(k,m)
    ~ 
    \mu^m\\
&< \infty.
\end{align*}

\end{proof}

In this special case of an infinite sequence of independent Bernoulli random variables, an expression for the moments involving only the first moments of the $Y_i$s and of $X$ can be easily had as well.
\begin{proposition}
\label{prop:infinite_moments}
	Let $(Y_i)_{i\ge 1}$ be a sequence of \textbf{independent} Bernoulli$(p_i)$ random variables with $\sum_{i\ge 1} p_i = \mu < \infty$. For the  Bernoulli (infinite) sum random variable $X = \sum_{i\ge 1} Y_i$, and for any $k \ge 2$, the $k$th \textbf{moment} of $X$ is
	\[E(X^k) = \sum_{m=1}^{k}S(k,m) \left[
	\mu^{m} - \sum_{s= 0}^{k-2}
		\mu^s(m-1-s)\sum_{\{i_1\,\ldots,i_{m-1-s}\}\subset \Naturals}
			p_{i_1}^2(p_{i_2}\dotsb p_{i_{m-1-s}})
	\right].\]
	\end{proposition}

\begin{proof}
	Fix an integer $r\geq 2$ and note that
\begin{align*}
\sum_{\{i_1, \ldots, i_{r}\}\subset \Naturals}p_{i_1}\dotsb p_{i_r}  
		&= \sum_{\{i_1, \ldots, i_{r-1}\}\subset \Naturals}
			p_{i_1}\dotsb p_{i_{r-1}}
			\left(
				\sum_{i \not\in \{i_1, \ldots, i_{r-1}\}}
					p_{i}  
						\right)\\
		&= \sum_{\{i_1, \ldots, i_{r-1}\}\subset \Naturals}
			p_{i_1}\dotsb p_{i_{r-1}}
			\left(
					\mu - \sum_{i \in \{i_1, \ldots, i_{r-1}\}} p_{i}  
						\right)\\
		&= \mu\sum_{\{i_1, \ldots, i_{r-1}\}\subset \Naturals}
				p_{i_1}\dotsb p_{i_{r-1}}
					- (r-2)\sum_{\{i_1, \ldots, i_{r-1}\}\subset \Naturals} 
						p_{i_1}^2\dotsb p_{i_{r-1}},					
\end{align*}
where the last equality follows from the fact that 
	$$\sum_{\{i_1, \ldots, i_{r-1}\}\subset \Naturals}  p_{i_1}^2 p_{i_2} \dotsb p_{i_{r-1}} = \sum_{\{i_1, \ldots, i_{r-1}\}\subset \Naturals}  p_{i_1} p_{i_2}^2 \dotsb p_{i_{r-1}} = \dotsb = \sum_{\{i_1, \ldots, i_{r-1}\}\subset \Naturals}  p_{i_1} p_{i_2}\dotsb p_{i_{r-1}}^2.$$
Recursively rewriting $\sum_{\{i_1, \ldots, i_{r}\}\subset \Naturals}p_{i_1}\dotsb p_{i_r}$ in terms of sums over one fewer index (viz.,  $r-1$ indices) each time gives the desired result via Proposition \ref{prop:bern_moments}.
\end{proof}

\subsection{Factorial moments}
The $k$th falling factorial of $x$ is the $k$th degree polynomial in $x$
\[
\ffact{x}{k} := x(x-1)(x-2)\dotsb (x-(k-1))= \prod_{m=0}^{k-1}(x-m),
\]
where $k \in \Naturals$ and $x \in \Reals$.
Replacing $x$ by a random variable $X$, the corresponding $k$th \textbf{factorial moment} is defined to be $E\left(\ffact{X}{k}\right)$.  Like $E(X^k)$ this is the expected value of the product of $k$ terms.  Note this is different from $E(X!)$, the \textbf{expected factorial} of $X$,  where the number of products in $X!$ is itself be a random variable (viz., $X$).

%

The $k$th power of $x$ can be expressed \cite[p. 74]{stanley2011enumerative} in terms of falling factorials as
\[x^k = \sum_{m = 1}^{k} S_2(k, m)\ffact{x}{m},\]
where $S_{2}(k,m)$ is the Stirling number of the second kind,
typically defined as the number of ways to partition a set of $k$ labelled objects into $m$ nonempty unlabelled subsets.
It follows, then, that these are directly related to the number of surjections from a $k-$set onto an $m-$set as 
\[S(k, m) = m! ~S_2(k, m)\]
and hence that
\[x^k = \sum_{m = 1}^{k} \frac{S(k, m)}{m!} \ffact{x}{m} = \sum_{m = 1}^{k} S(k, m) {x \choose m}.\]
Similarly, the falling factorial is written as a sum of powers as
\[[x]_k = \sum_{m = 1}^k S_1(k, m) x^m\]
where $S_1(k,m)$ is the Stirling number of the first kind.
Similar expressions may now be found involving a Bernoulli sum $X$ in place of $x$.

First, we relate $X \choose m$ to the Bernoullis that define $X$.
\begin{proposition}
	\label{prop:bern_choose_function}
	If $X$ is a Bernoulli sum $X = \sum_{i\in I}Y_i$ for some countable indexing set $\set{I}$, then for $m\geq 1$,
	\[\binom{X}{m}
	= 
	\sum_{\{i_1,\ldots, i_m\}\subseteq ~\set{I}} 
		Y_{i_1} \dotsb Y_{i_m}.\]
\end{proposition}

\begin{proof}
Let $\set{J} = \{i \in \set{I} : Y_i=1\}$ denote the subset of the indices in $\set{I}$ for which $Y_i = 1$. Then
\begin{align*}
\sum_{\{i_1,\ldots, i_m\}\subseteq \set{I}} 
		Y_{i_1} \dotsb Y_{i_m} &= 
			\sum_{\{i_1,\ldots, i_m\}\subseteq ~\set{J}} 
		Y_{i_1} \dotsb Y_{i_m}\\
	  &= \sum_{\{i_1,\ldots, i_m\}\subseteq ~\set{J}} 
		1\\
	  &= \binom{|\set{J}|}{m}\\
	  &= \binom{X}{m}.
\end{align*}
 	
\end{proof}
\noindent
An earlier, inductive, proof of this result for the case of finite $\set{I}$ is given by \citeA{iyer1958theorem}.

A similar approach yields a general result relating $X!$ to its Bernoulli constituents.
\begin{proposition}
	\label{prop:expected_factorial} Let $X = \sum_{i\in \set{I}} Y_i$ be a Bernoulli sum, where $\set{I}$ is a countable indexing set. Then we may write $X!$ in terms of the $(Y_i)$ as follows
	$$X! = \sum_{\set{H} \subseteq \set{I}} |\set{H}|!  ~
				\left( \prod_{i\in \set{H}}
					Y_{i} ~
						\prod_{i\not\in \set{H}}
							(1-Y_i) \right).$$

\end{proposition}
\begin{proof}
	Consider the set $\set{J} := \{i \in \set{I}: Y_i = 1\}$. In this case, $|\set{J}| = X$ and  $|\set{J}|! = X!$.
       For any other set $\set{H} \subseteq \set{I}$, either $\set{H} = \set{J}$,  or $\set{H} \ne  \set{J}$.
	
	If   $\set{H} = \set{J}$, then
	\[\prod_{i \in H} Y_{i} ~
	  \prod_{i \not\in \set{H}} (1-Y_i) 
	  =  \prod_{i \in \set{H}} 1 ~ \prod_{i \not\in \set{H}}(1-0) 
	  = 1 \]
	  
	If  $\set{H} \ne  \set{J}$, then there exists $j$ for which $ j \in \set{J} $ but $j \not\in \set{H}$. Then,
	\[
	 \prod_{i \in \set{H}} Y_{i} ~
	 \prod_{i \not\in \set{H}}(1-Y_i) = \prod_{i\in \set{H}}Y_i \times 0 = 0.
	 \]
	 
	 Together these give
	 \[\sum_{\set{H} \subseteq \set{I}} |\set{H}|!  ~
				\left( \prod_{i\in \set{H}} Y_{i} ~
				\prod_{i\not\in \set{H}} (1-Y_i) \right)= |\set{J}|! \prod_{i\in \set{J}} Y_{i} ~
				\prod_{i\not\in \set{J}} (1-Y_i) = |\set{J}|! = X!
				\]
\end{proof}

Taking expectations yields the following expressions for a Bernoulli sum $X = \sum_{i \in \set{I}} Y_i$:
\begin{itemize}
\item the $k$th factorial moment in terms of the Bernoulli random variables
\[ E(\ffact{X}{k}) = k!  \sum_{\{i_1,\ldots, i_k\} \subseteq ~ \set{I}} 
		E(Y_{i_1} \dotsb Y_{i_k})\]
		or, in terms of the moments of $X$ as
\[ E(\ffact{X}{k}) =  \sum_{m = 1}^k S_1(k, m) E(X^m) \]		
\item the $k$th moment in terms of the factorial moments of $X$
\[E(X^k) 
= \sum_{m=1}^{k} S_2(k, m)~ E(\ffact{X}{m}) \]
\item the $k$th central moment 
\[ E((X - \mu)^k) 
    = (-\mu)^{k} +   \sum_{j=1}^{k}  \left( \sum_{m=j}^{k} {S_2(m, j)} \binom{k}{m} (-\mu)^{k-m}  \right) E(\ffact{X}{j})
\]

\item and the expected factorial in terms of the Bernoulli random variables
\[  E(X!) = \sum_{\set{H}\subseteq \set{I}} |\set{H}|!   ~ \times
					E\left[\prod_{i\in \set{H}} 
						Y_{i} ~
							\prod_{i\not\in \set{H}}
								(1-Y_i)\right]. \]
\end{itemize}
Central moments for small $k$ can always be written in terms of the moments or in terms of the factorial moments.  When $k =2$, a nice symmetry appears in either expression for the variance of $X$:
\[ Var(X) = E(X^2) - (E(X))^2 = E([X]_2) - [E(X)]_2.\]

\subsection{A statistical interpretation}
Central moments are statistically meaningful for any random variable $X$ where available.  However, when $X$ is a Bernoulli sum a few more meaningful interpretations are available.

Imagine a collection of individuals $i \in \set{I}$, from which a random number $X$ provides a population $\set{J}$ of size $X$.  Samples of fixed size $k$ are to be drawn from the resulting population $\set{J}$.    Here,  $X = \sum_{i \in \set{I}} Y_i$ and $\set{J}  = \{i \in \set{I}: Y_i = 1\}$ with (possibly dependent) random variables $Y_i \follows Bernoulli(p_i)$ (indicating inclusion in the population $\set{J}$ when $Y_i =1$ and exclusion when $Y_i = 0$).

In this case, the \textit{expected number of samples of size} $k$ 
\begin{itemize}
\item is the $k$th factorial moment $E(\ffact{X}{k})$ when sampling \textit{without replacement} and
\item is the $k$th moment $E(X^k)$ when sampling \textit{with replacement}.
\end{itemize}
The expected factorial $E(X!)$ is the \textit{expected number of permutations} one would have in the indices found by forming a population in this way.

\subsection{Generating functions}

Various generating functions for a Bernoulli sum $X$ are now easily had by substitution of

\begin{itemize}
\item $E(X^k)$ in the \textit{moment generating function}
 \[ M_X(s) = E(e^{sX}) = 1 + \sum_{k=1}^{\infty}E(X^k)\frac{s^{k}}{k!} \]
\item $E(\ffact{X}{k})$ in the \textit{factorial moment generating function} \cite<e.g., see>[p. 59]{johnson2005univariate}
\[H_X(s) = 1 + \sum_{k=1}^{\infty}E(\ffact{X}{k})\frac{s^{k}}{k!} \]
\item and, from \citeA{frechet1943extension},
\[
 Pr(X = x) = \sum_{j\geq x} (-1)^{x+j}\binom{j}{x}\frac{E(\ffact{X}{j})}{j!} ,
\]
or, after substitution for the factorial moments, 
\[Pr(X = x) = \sum_{j\geq x} (-1)^{x+j}\binom{j}{x}\sum_{\{i_1,\ldots, i_j\} \subseteq ~ \set{I}} 
		E(Y_{i_1} \dotsb Y_{i_j}),\]
 the probability $Pr(X = x)$ into the \textit{probability generating function}
\[G_X(s) = E(s^X) = \sum_{k=0}^{\infty}s^k Pr(X=k).\]
\end{itemize}
The factorial moment generating function, $H_X(s)$, can be related 
\cite<again, see>[p. 59]{johnson2005univariate} to the probability generating function, $G_X(s)$,  as
\begin{align*}
 H_X(s) = G_X(1+s) = E((1+s)^X).
\end{align*}
It follows that whenever factorial moments are such that $H_X(s)$ has a tidy closed form, the probability generating function of $X$ might be easily obtained through the reverse relation
\begin{equation}
G_X(s) = H_X(s-1).
\label{eq:fmgf}
\end{equation}
This approach will be illustrated for the binomial distribution in Section \ref{sec:binomial}, and  
for the classic matching problem of Section \ref{sec:matching_problem},  to determine expressions for the probability generating function of $X$ in each of these classic cases.

\section{Classic examples} 
\label{sec:classics}

Bernoulli sums naturally arise in many classic problems and lead to well known distributions.
In this section, the results of Section \ref{sec:moments} are applied to several of these where the Bernoulli sum is over a finite index set (of size $n$), namely
\[ X = \sum_{i = 1}^{n} Y_i \]
where $Y_i \follows Bernoulli(p_i)$ with $p_i = Pr(Y_i =1) = 1 - Pr(Y_i = 0)$ for $i = 1, \ldots, n$.

\subsection{Binomial $X$}
\label{sec:binomial}

The simplest case where the $Y_i$s are independent and identically distributed (i.i.d.) with $p_i = p ~ \forall ~i$, $X \follows binomial(n,p)$. 
The $k$th moment of $X$ can be written as 
\begin{align*}
\label{binom_moments}
E(X^k) 
& = \sum_{m=1}^{k}
	 S(k,m)
	 \sum_{\{i_1, \ldots,  i_m \}\subseteq \Nset{n}}
	~ p^m\\
& = \sum_{m=1}^{k}
	 S(k,m)
	\binom{n}{m}
	~ p^m . \\
	\end{align*}
\citeA{knoblauch2008closed} found an equivalent expression through a recursive argument. 
Here the result is easily had from the simple application of the more general Proposition \ref{prop:bern_moments}.
Central moments follow 
from Proposition \ref{prop:bern_central_moments}:
\begin{align*}
E((X-\mu)^{k}) 
&= (-np)^k + 
	\sum_{\ell=1}^{k}
	(-np)^{k-\ell}
	\sum_{m=1}^{\ell}
	S(\ell,m)
	\binom{n}{m}
	p^m.
\end{align*}
The $k$th factorial moment has a appealingly simple expression $E(\ffact{X}{k}) = [n]_k ~ p^k$ derived as
\[
E(\ffact{X}{k}) = 
					k! \sum_{\{i_1,\ldots, i_k\}\subseteq \Nset{n}}
							p_{i_1} \dotsb p_{i_k}\\
							= 
					k! \binom{n}{k} p^k = [n]_k ~ p^k,
\]
the familiar $E(X) = np$ being the special case when $k = 1$. 

Note that whenever $k \ge n$, many terms disappear in the above moment expressions since 
$S(n, m)$ vanishes whenever  $m > n$ and the sum $\sum_{\{i_1,\ldots, i_k\}\subseteq \Nset{n}
}$ is over the empty set. 

The moment generating function of a binomial $X$ also has a new expression following application of Equation \ref{eq:mgf},
namely
\begin{equation}
M_X(t) = 1 + \sum_{k \ge 1} \frac{t^k}{k!} 
   ~ \sum_{m=1}^{k}
	 ~S(k,m) 
	 \binom{n}{m}
	 p^m
\label{eq:mgf}
\end{equation}
compared to $M_X(t) = (1 - p + p e^t)^n$.

Recall that the probability generating function of $X \follows binomial(n,p)$ is 
\[G_X(s) = \left((1-p) + ps \right)^n. \]
By Equation \ref{eq:fmgf}, we find that the factorial moment generating function for $X$ is
\begin{align*}
  H_X(s) &= ((1-p) + p(1+s))^n\\
  &= \sum_{m=0}^{n} s^m
  	\left(\sum_{r=m}^{n}
  		  \sum_{\ell = 0}^{n-r} 
  			\binom{n}{r}
  				\binom{r}{m}
  		 			(-1)^{\ell}p^{\ell + r}
  						\right),
\end{align*}
by applying the binomial theorem and changing the order of summation. This provides us an additional expression for the $k-$th factorial moment of $X$:
\begin{align*}
 \ffact{n}{k}p^k = \sum_{r=k}^{n}
  		  \sum_{\ell = 0}^{n-r} 
  			\binom{n}{r}
  				\binom{r}{k}
  		 			(-1)^{\ell}p^{\ell + r}.
\end{align*}

\subsubsection{Poisson binomial $X$}
If $Y_i \follows  Bernoulli(p_i)$ independently for all $i$ but $p_i \ne p_j$ for (at least one) $i \ne j$, the distribution of $X$ is called a Poisson binomial 
distribution \cite<e.g., see>{shah1973distribution}.
%
The various moments of $X$ are exactly as given by the relevant results of Section \ref{sec:moments} with $E(Y_{i_1} \dotsb Y_{i_m})$ everywhere replaced by $p_{i_1} \dotsb p_{i_m}$. 
So too for its moment generating function.

\subsection{Hypergeometric $X$}

Consider a sample of size $n$ randomly drawn without replacement from a population of $N$ individuals where $g$ of them have some trait which is absent from the remaining $N - g$.
The indicator random variable, $Y_i$, records 
 if the $i$th individual selected has the desired trait ($Y_i =1$) or not ($Y_i =0$) and $X = \sum_{i=1}^{n} Y_i$ counts the number in the sample having the trait.
 
The $i$th draw will be a Bernoulli random variable $Y_i$ with probability 
\[p_i = \frac{g - \ell}{N - (i-1)}\]
where $\ell$ is the number of previous $(i-1)$ draws having the trait.
A sample of $m$ of these Bernoullis will have
\[E(Y_{i_1}\dotsb Y_{i_m}) 
     = \frac{g(g-1)\dotsb (g-m+1)}
               {N(N-1)\dotsb (N-m+1)}
                = \frac{\ffact{g}{m}}
                          {\ffact{N}{m}}
               \]
provided $m \le g$ and will be zero whenever $m > g$ (since at least one $Y_i$ must be zero).

The $k$th moment of $X$, following Proposition \ref{prop:bern_moments}, is now
\begin{align*}
  E(X^k) &= 
	\sum_{m=1}^{\min{k, g}}S(k,m)
		\sum_{\{i_1, \ldots,  i_m \}\subseteq \Nset{n}}
			\frac{\ffact{g}{m}}
                          {\ffact{N}{m}}\\
	&= \sum_{m=1}^{\min{k, g}}
			S(k,m)
				\binom{n}{m} 
					\frac{\ffact{g}{m}}
                          {\ffact{N}{m}},
\end{align*}
where the last equality followed from $m-$symmetry.
%
Similarly, the central moments are 
\begin{align*}
	E((X-\mu)^k) 
	&= \left(-n\frac{g}{N}\right)^k 
	   + \sum_{\ell=1}^{M} 
	       \binom{k}{\ell}\sum_{m=1}^{\ell}
	        S(\ell,m)\binom{n}{m}\frac{\ffact{g}{m}}
                          {\ffact{N}{m}}
                          \left(- n\frac{g}{N}\right)^{k-\ell}
\end{align*}
where $M = \min{k, g}$.

The factorial moments again have a pleasingly simple expression when $k \le g$ (zero whenever $k > g$), namely,
\begin{align*}
  E(\ffact{X}{k}) 
  &= ~k!
  		\binom{n}{k}
			\frac{\ffact{g}{k}}{\ffact{N}{k}}
  ~=~ \ffact{n}{k} ~\frac{\ffact{g}{k}}{\ffact{N}{k}}.
\end{align*}
Where the binomial $E(\ffact{X}{k}) = \ffact{n}{k}p^k$, the hypergeometric now has $\frac{\ffact{g}{k}}{\ffact{N}{k}}$ in place of $p^k$, as one might expect.

\subsection{CMP-binomial $X$}
\label{sec:CMP_binomial}

For $n \in \Naturals, p \in [0, 1], \nu \in \Reals$, a random variable $X$ has a \textit{Conway-Maxwell-Poisson (CMP) binomial} distribution with parameters $(n, p, \nu)$ if its probability mass at $X = j$ ($j \in \Nset{n}$) is given by
\begin{align*}
  Pr(X=j) &= \frac{1}{C_{n,p,\nu}}\binom{n}{j}^\nu p^j(1-p)^{n-j},
\end{align*}
where $C_{n,p,\nu}$ is the normalizing constant
\begin{align*}
  C_{n,p,\nu} &= \sum_{j = 0}^{n}\binom{n}{j}^\nu p^j(1-p)^{n-j}.
\end{align*}
The distribution is formed from a Conway-Maxwell-Poisson (CMP) random variable conditional on the sum of that variable and another one independently generated from a different CMP-distribution. 

Just as the CMP-distribution generalizes a Poisson random variable to model count data having variability larger (over dispersed) or smaller (under dispersed) than that of a Poisson, the \textit{CMP-binomial} generalizes the binomial distribution.  A CMP-binomial distribution is binomial when $\nu =1$ and has larger (smaller) variance than a binomial when $\nu < 1$ ($\nu > 1$).  When $\nu = 0$, the most extreme values of $0$ and $n$ are favoured;  when  $\nu \to \infty$ the count $X$ achieves the middle value of $n/2$ when $n$ is even and  $(n \pm 1)/2$ when $n$ is odd.  See \citeA{shmueli2005useful} for details.

\citeA{shmueli2005useful} remark that the random variable $X$ can also be viewed as a sum of exchangeable, Bernoulli random variables $Y_i$ with joint probability 
\begin{align*}
  Pr(Y_1 = y_1,\ldots, Y_n = y_n) = \frac{1}{C_{n, p,\nu}} \binom{n}{\sum_{i=1}^{n}y_i}^{\nu-1}p^{\sum_{i=1}^{n}y_i}(1-p)^{n-\sum_{i=1}^{n}y_i},
\end{align*}
where $\nu > 1$ in the case of negatively correlated trials and $\nu < 1$ for positively correlated trials. This observation allows an expression to be written for the moments of $X$ from an expression $Pr(Y_{i_1} = 1, \ldots, Y_{i_m}=1)$ for an arbitrary $m-$set $\{i_1,\ldots, i_m\}\subseteq \Nset{n}$.
\begin{align*}
  Pr(Y_{i_1},\ldots, Y_{i_m}) &= 
  \sum_{\substack{y_j \in \{0,1\}\\ \forall j \not\in \{i_1,\ldots, i_m\}}} 
  	Pr(Y_{i_1}=1,\ldots, Y_{i_m}=1, \text{ and } Y_j = y_j, \forall j \not\in \{i_1, \ldots, i_m\})\\
  	&= 
  \sum_{\substack{y_j \in \{0,1\}\\ \forall j \not\in \{i_1,\ldots, i_m\}}} 
  	\frac{1}{C_{n, p, \nu}}\binom{n}{m+\sum_{j\not\in\{i_1,\ldots, i_m\}y_j}}^{\nu-1}
  		\\
  		&~~~~~~~~~~~~~~~~~~~~~~~~~~~~~~~~~~~~~~\times p^{m+\sum_{j\not\in\{i_1,\ldots, i_m\}}y_j}  (1-p)^{n-(m+\sum_{j\not\in\{i_1,\ldots, i_m\}}y_j)}\\
		&\\
  	&= \frac{1}{C_{n, p, \nu}} 
	\sum_{s = 0}^{n-m} 
  	\sum_{\substack{y_j \in \{0,1\}\\ \forall j \not\in \{i_1,\ldots, i_m\}\\ \sum_{j\not\in\{i_1,\ldots, i_m\}}y_j = s}} 
  		\binom{n}{m+s}^{\nu-1}
  		\times p^{m+s}(1-p)^{n-(m+s)}\\
		&\\
  	&= \frac{1}{C_{n, p, \nu}}
  		\sum_{s = 0}^{n-m} 
  			\binom{n-m}{s}
  			\binom{n}{m+s}^{\nu-1}
  				p^{m+s}(1-p)^{n-(m+s)}\\
				&\\
  	&= \frac{1}{C_{n, p, \nu}}
  		\sum_{\ell = m}^{n} 
  			\binom{n-m}{\ell-m}
  			\binom{n}{\ell}^{\nu-1}
  				p^{\ell}(1-p)^{n-\ell}.\\  				
\end{align*}
The $k$th moment of $X \follows CMP-binomial(n, p, \nu)$  is then
\begin{align*}
  E(X^k) 
  	&= 
  		\sum_{m=1}^{k} S(k,m) 
  			\sum_{\{i_1,\ldots, i_m\}\subseteq \Nset{n}}
  				E(Y_{i_1}\dotsb Y_{i_m})\\
  	&=  \frac{1}{C_{n,p,\nu}}
	      \sum_{m=1}^{k} S(k,m) 
  			\binom{n}{m}
  		\sum_{\ell = m}^{\min{n, k}} 
  			\binom{n-m}{\ell-m}
  			\binom{n}{\ell}^{\nu-1}
  				p^{\ell}(1-p)^{n-\ell}\\
  	&=  \frac{1}{C_{n,p,\nu}}
	      \sum_{m=1}^{k} S(k,m) 
  		\sum_{\ell = m}^{\min{n, k}}
  			\binom{\ell}{m}
  			\binom{n}{\ell}^{\nu}
  				p^{\ell}(1-p)^{n-\ell}.\\  
\end{align*}
Similarly, the $k$th central moment is
\begin{align*}
\left(-np \right)^{k} +
         \frac{1}{C_{n,p,\nu}}\sum_{\ell=1}^{k}
       \binom{k}{\ell} 
         \left(-np\right)^{k-\ell}
        \sum_{m=1}^{\ell} 
         S(\ell,m)
  					\sum_{\ell = m}^{\min{n, k}} 
  						\binom{\ell}{m}
  							\binom{n}{\ell}^{\nu}
  								p^{\ell}(1-p)^{n-\ell}
\end{align*}
and the $k$th factorial moment
\begin{align*}
E(\ffact{X}{k}) =
\frac{k!}{C_{n,p,\nu}} 
  		\sum_{\ell = k}^{n} 
  			\binom{\ell}{k}
  			\binom{n}{\ell}^{\nu}
  				p^{\ell}(1-p)^{n-\ell} .
\end{align*}

\subsection{The empty urns problem}
\label{sec:empty_urns}

Consider the problem of assigning $\ell$ indistinguishable balls uniformly at random into $n$ distinguishable urns. 
Let $Y_i$ be 1 if urn $i$ is empty and 0 otherwise, and  $X = \sum_{i=1}^{n}Y_i$ be the Bernoulli sum 
counting the total number of empty urns. 

Through a straightforward counting argument, it can be shown that there are $\binom{\ell + n -1}{n}$ ways to distribute $\ell$ indistinguishable balls into $n$ distinguishable urns and therefore
\begin{align*}
  Pr(Y_i = 1) &= \frac{\# \text{ ways to distribute $m$ balls into $n-1$ urns}}{\# \text{ ways to distribute $m$ balls into $n$ urns}}\\
  &= \frac{\binom{n+\ell-2}{\ell}}{\binom{n+\ell-1}{\ell}}\\
  &= \frac{n-1}{\ell+n-1}.
\end{align*}
By the same argument, for a subset $\{i_1, \ldots, i_m\}$ of $\Nset{n}$,
\begin{align*}
  Pr(Y_{i_1} = 1, Y_{i_2} = 1, \ldots, Y_{i_m}=1)
  &= \frac{(n-1)(n-2)\dotsb (n-m)}{(\ell+n-1)(\ell+n-2)\dotsb (\ell+n-m)} 
     = \frac{\ffact{n-1}{m}}{\ffact{\ell + n -1}{m}}.
\end{align*}
The $k$th moment of $X$ is 
\begin{align*}
  E(X^k) &= 
	\sum_{m=1}^{k}\sum_{\{i_1, \ldots,  i_m \}\subseteq \Nset{n}} S(k,m)E(Y_{i_1}\dotsb Y_{i_m})\\
	&= \sum_{m=1}^{k}S(k,m)\binom{n}{m}\frac{\ffact{n-1}{m}}{\ffact{\ell + n -1}{m}},
\end{align*}
the $k$th central moment  (from Proposition \ref{prop:bern_central_moments})
\begin{align*}
  E((X-\mu)^k) 
&= 
	\left(-\mu\right)^{k} +
\sum_{\ell=1}^{k}
       \binom{k}{\ell} 
         \left(-\mu\right)^{k-\ell}
        \sum_{m=1}^{\ell} 
         S(\ell,m)
        \binom{n}{m} 
        	\frac{\ffact{n-1}{m}}{\ffact{\ell+n-1}{m}},
\end{align*}
where $\mu = n\frac{n-1}{\ell+n-1}$.  
The $k$th factorial moment has a particularly simple representation as
\begin{align*}
  E(\ffact{X}{k}) 
  	&=
  		\frac{\ffact{n}{k}\ffact{n-1}{k}}{\ffact{\ell+n-1}{k}}.
\end{align*}
Matching moments shows $X$ of the urn problem to have a Hypergeometric distribution with parameters $(n, n-1, \ell + n - 1)$. 

%
\subsection{The matching problem}
\label{sec:matching_problem}
The matching problem dates back to Montmort \cite{montmort1713essay} and is the problem of taking $n$ paired elements, randomly permuting the first elements over all pairs, then letting $X$ be the number of correctly matched pairs after the random permutation.  Examples are $n$ letters matched correctly to $n$ envelopes, couples separated at a dance and dance partners formed by randomly assigning one of each sex to the pair, and so on.
The random variable $X$ can be expressed as a sum of Bernoulli random variables taking value 1 when a correct match occurs and zero otherwise.

More abstractly,  let  
%
$f:\Nset{n}\to \Nset{n}$ be a permutation on $\Nset{n}$ and let $X$ denote the number of fixed points of $f$ (i.e.,  the number of $i\in \Nset{n}$ for which $f(i) = i$). 
If $f$ is picked uniformly at random from the set of all permutations on $\Nset{n}$, denoted $Sym(n)$, then the distribution of $X$ can be shown to tend to $Poisson(1)$ as $n \to \infty$. 
We do that by expressing $X$ as a sum of Bernoullis and then examining and comparing moments.

Let $Y_i$ be the Bernoulli random variable recording if $f(i) = i$. As $f$ is chosen uniformly at random from $Sym(n)$, 
$$Pr(Y_i = 1) = \frac{1}{n}$$ 
for all $i = 1,\ldots, n$. 

Fix $m\leq n$ and consider an $m-$subset $\{i_1,\ldots, i_m\} \subseteq \Nset{n}$. The probability that all $\{i_1,\ldots, i_m\}$ are fixed points of $f$ is $(n-m)!/n!$. This is because if $f(i_j) = i_j$ for all $j\in\{1,\ldots, m\}$, one must only consider how to assign the other $(n-m)$ points so that $f$ is a bijection. There are $(n-m)!$ ways to do this as $|Sym(n-m)| = (n-m)!$. Therefore, 
$$Pr(Y_{i_1} = 1,\ldots, Y_{i_m} =1) = \frac{(n-m)!}{n!}.$$
By Proposition $\ref{prop:bern_moments}$, for $k\leq n$
\begin{align*}
  E(X^k) &= \sum_{m=1}^{k}
  				S(k,m)
  					\sum_{\{i_1,\ldots, i_m\}\subseteq \Nset{n}}
  						Pr(Y_{i_1} = 1,\ldots, Y_{i_m} = 1)\\
  		 &=	\sum_{m=1}^{k}
  				S(k,m)
  					\sum_{\{i_1,\ldots, i_m\}\subseteq \Nset{n}}
  						\frac{(n-m)!}{n!}\\
  		 &=	\sum_{m=1}^{k}
  				S(k,m)
  					\binom{n}{m}
  						\frac{(n-m)!}{n!}\\
  		 &=	\sum_{m=1}^{k}
  				\frac{S(k,m)}{m!}\\
  		 &= \sum_{m=1}^{k}S_2(k,m)\\
  		 &= B_k,				  						
\end{align*}
where $B_k$ is the $k-$th \textit{Bell number}, the number of ways to partition a set of size $k$ into a family of nonempty, unlabelled, pairwise disjoint subsets. On the other hand, for $k> n$, the inner sum vanishes whenever $m > n$ and hence  
\begin{align*}
    E(X^k) &= \sum_{m=1}^{k}
  				S(k,m)
  					\sum_{\{i_1,\ldots, i_m\}\subseteq \Nset{n}}
  						Pr(Y_{i_1} = 1,\ldots, Y_{i_m} = 1)\\
  		 &= \sum_{m=1}^{n}
  				S(k,m)
  					\sum_{\{i_1,\ldots, i_m\}\subseteq \Nset{n}}
  						Pr(Y_{i_1} = 1,\ldots, Y_{i_m} = 1)\\
  		 &=	\sum_{m=1}^{n}
  				\frac{S(k,m)}{m!}\\
  		 &= \sum_{m=1}^{n}S_2(k,m)\\
  		 &= B_k - \sum_{m=n+1}^{k}S_2(k,m),		
\end{align*}
which can be interpreted as the number of ways to partition a set of size $k$ into at most $n$ classes.

For $k\leq n$, the $k-$th factorial moments of $X$ is given by
\begin{align*}
  E(\ffact{X}{k}) &= k! \sum_{\{i_1,\ldots, i_k\}\subseteq \Nset{n}}
  						Pr(Y_{i_1} = 1,\ldots, Y_{i_k} = 1)\\
  				  &= k! \frac{(n-k)!}{n!}\binom{n}{k}\\
  				  &= 1.
\end{align*}
Therefore, the factorial moment generating function of $X$ is given by
\[H_X(s) = \sum_{k=0}^{n}\frac{s^k}{k!},\]
which by Equation \ref{eq:fmgf} gives us that the probability generating function of $X$ is 
\begin{align*}
G_X(s) &= H_X(s-1)\\ 
&= \sum_{k=0}^{n}\frac{(s-1)^k}{k!}\\
&= \sum_{\ell = 0}^{n}s^{\ell}\sum_{k = \ell}^{n} \binom{k}{\ell}\frac{(-1)^{k-\ell}}{k!},
\end{align*}
which provides us with another method for deriving the probability distribution of $X$ \textit{without} using Principle of Inclusion/Exclusion or counting derangements in permutations.

Now,  the exponential generating function for $B_k$ given by  \cite<e.g., see>[p. 74]{stanley2011enumerative} 
$$\sum_{k\geq 0} B_k \frac{t^k}{k!} = e^{e^t-1}$$
is a special case (viz., $\lambda = 1$) of the moment generating function for a Poisson($\lambda$) random variable $W$ 
$$M_W(t) = \sum_{k\geq 0} E(W^k)\frac{t^k}{k!} =  e^{\lambda(e^t-1)}.$$
For $k\leq n$, the moments of $X$ match those of $W \follows$ Poisson(1);  as $n \to \infty$, the moment generating functions agree and $X$ converges to a Poisson(1) random variable.

\subsection{The Poisson limit of a binomial}
\label{sec:binom_to_poisson}
Proposition $\ref{prop:bern_moments}$ can also be used to provide a novel proof that a binomial$(n, p)$ random variable approaches a Poisson($\lambda$) with $np \to \lambda$  as $n \to \infty$. 

For $X = \sum_{i=1}^{n}Y_i$,  each $Y_i$ is independent, identically distributed Bernoulli$(p)$ random variable and, as seen earlier, the $k$th moment
\begin{align*}
E(X^k) &= \sum_{m=1}^{k}
	 S(k,m)
	\binom{n}{m}
	~ p^m\\
	&= \sum_{m=1}^{k}
	 S(k,m)
	\binom{n}{m}
	~ \left(\frac{np}{n}\right)^m\\
	&= \sum_{m=1}^{k}
	 \frac{S(k,m)}{m!}
	\left(\frac{\ffact{n}{m}}{n^m}\right)\left(np\right)^m.
\end{align*}
As $n \to \infty$, the ratio $\left(\frac{\ffact{n}{m}}{n^m}\right) \to 1$, $np \to \lambda$,  and 
\[E(X^k) \to \sum_{m=1}^{k}
	 \frac{S(k,m)}{m!}
	\lambda^m = \sum_{m=1}^{k}
	 S_2(k,m)
	\lambda^m \]
which is the $k$th moment of a Poisson($\lambda$) random variable expressed as a Touchard polynomial in $\lambda$ \cite<e.g., see>{riordan1937moment}.  It follows that as $n \to \infty$, $X \follows binomial(n, p)$ converges to a Poisson($\lambda$) with $\lambda = \lim_{n \to \infty} np$.

When $p = \frac{1}{n}$, then  binomial $X$ converges to Poisson(1) and its $k$th moment is the Bell number $B_k$, as in the matching problem of Section \ref{sec:matching_problem}.

%

\section{Counts more generally}
\label{sec:counts_generally}
Focus has been on Bernoulli sums that rise naturally as counts of events.
In this section, we consider any ``count random variable" $N$ to be that having support on \textit{any} subset of the extended natural numbers $\NaturalsZero$ (e.g., counts $n$ would be impossible whenever $Pr(N=n) = 0$).  Previous results are extended to $N$ by matching it to a Bernoulli sum $X$ constructed from the upper tail probabilities of $N$. 
%


\begin{proposition}
\label{prop:discrete_rvs}
	Let $N$ denote any discrete random variable on $\NaturalsZero$. Consider the Bernoulli random variable $Y_i$ indicating whether $N \geq i$ or not; that is 
	\[ Y_i = \left\{ 
	  \begin{array}{lcl}
	  1 & \text{if}  & N\geq i ,\\
	  &&\\
	  0 && \text{otherwise.}
	  \end{array}
	  \right. \]
Then the Bernoulli sum $X = \sum_{i = 1}^{\infty} Y_i$ has the same distribution as $N$.
\end{proposition}

\begin{proof}
\begin{align*}
Pr(X = x) &= 
	Pr\left(\sum_{i=1}^{\infty} Y_i = x\right)\\ 
		  &= 
		  	Pr([Y_i = 1, \forall i: i \leq x]\cap [Y_i = 0,\forall i: i > x])\\ 
		  &= Pr(N = x).
\end{align*}

\end{proof}
The moments of $N$ are identified with those of $X$ and the Bernoullis $Y_i$ defined above.  The results for the moments and factorial moments now follow.
\begin{proposition}
\label{prop:discrete_moments}
	For any discrete random variable $N$ on $\NaturalsZero$, the $k$th moment of $N$ is
\begin{align*}
E(N^k) &= \sum_{m = 1}^{k}
			S(k,m) \sum_{M\geq m}\binom{M-1}{m-1}Pr(N \geq M),
\end{align*}
and the $k$th factorial moments is
\begin{align*}
 E(\ffact{N}{k}) &= k! 
 					\sum_{M\geq k}\binom{M-1}{k-1}Pr(N \geq M).
\end{align*}

\end{proposition}
\begin{proof}
From Proposition \ref{prop:bern_moments}, 
\begin{align*}
  E(N^k) &= \sum_{m=1}^{k} 
  				S(k,m)
  					\sum_{\{i_1,\ldots, i_m\}\subset \Naturals} 
  						Pr(Y_{i_1} = 1, \ldots, Y_{i_m} = 1)
\end{align*}
Now,
\begin{align*}
  \sum_{\{i_1,\ldots, i_m\}\subset \Naturals} 
  						Pr(Y_{i_1} = 1, \ldots, Y_{i_m} = 1) &= 
	  \sum_{\{i_1,\ldots, i_m\}\subset \Naturals} 
	  	Pr(\max({i_1,\ldots, i_m}) \leq N)\\
	  	&= \sum_{M = 1}^\infty  \sum_{\{i_1,\ldots, i_{m-1}\}\subseteq M-1}
	  		Pr(N \geq M)\\
	  	&= \sum_{M=1}^{\infty}\binom{M-1}{m-1}Pr(N \geq M).	  			   						
\end{align*}
Therefore,
\begin{align*}
  E(N^k) &= \sum_{m=1}^{k} 
  				S(k,m)
  					 \sum_{M=1}^{\infty}\binom{M-1}{m-1}Pr(N \geq M).	  			   						
\end{align*}

For the factorial moment, 
\begin{align*}
  m!\sum_{M=m}^{\infty}\binom{M-1}{m-1}Pr(N \geq M) 
  	&= 
  		m!\sum_{M\geq m}\sum_{\ell \geq M} Pr(N = \ell) \binom{M-1}{m-1}\\
  	&= 
  		m!\sum_{\ell\geq m}Pr(N = \ell)
  			\sum_{M =1}^{\ell} 
  				\binom{M-1}{m-1}\\
  	&= 
  		m!\sum_{\ell\geq m}Pr(N = \ell)
  			\binom{\ell}{m}\\ 
  	&= 
  		m!\sum_{\ell\geq m}\frac{\ffact{\ell}{m}}{m!}Pr(N = \ell)\\
  	&= 
  		\sum_{\ell\geq m}\ffact{\ell}{m}Pr(N = \ell)\\
  	&= 	
  		\sum_{\ell\geq 0}\ffact{\ell}{m}Pr(N = \ell)\\  
  	&= E(\ffact{N}{m}) ~~~ \text{by definition}.					 				
\end{align*}
\end{proof}
Note that an expression for $E(N^k)$ has also recently been derived by \citeA[eq. (10)]{Chakra2019}, namely
\[E(N^k) = \sum_{i=0}^\infty ((i + 1)^k- i^k) Pr(N > i). \]
\citeA{Chakra2019} claim their formulation to be the first for $E(N^k)$ expressed in terms of the upper tail probability of $N$; if so, then Proposition \ref{prop:discrete_moments} may provide the second for $E(N^k)$ and the first for $E(\ffact{N}{k})$.
These results are 
best appreciated whenever the cumulative distribution function of $N$ has form allowing simplification, especially when multiplied by binomial coefficients.   

The remainder of this section 
explores application of Proposition \ref{prop:discrete_moments} to several  familiar cases.

\subsection{Geometric distribution}
\label{sec:geometric}
Suppose interest lay in the number $X$ of tosses of a coin at which the first head occurs; $X$ is a geometric($p$) distribution with $p$ being the probability of heads ($p = 0.5$ for a fair coin).  Surprisingly, the random variable $X$ can be written as a Bernoulli sum  $X =\sum_{i=1}^\infty Y_i$ for suitably defined Bernoulli $Y_i$s, allowing the previous results to be applied.

The representation is as follows.  Take $(Z_{i})_{i\geq1}$ to be the sequence of independent Bernoulli($p$) random variables representing the sequence of potential coin tosses ($Z_i = 1$ for heads and zero otherwise).  Let $N$ be the index of the first $Z_i = 1$ in the sequence, that is
\[N = \min{i: Z_i = 1}.\]
Consider now the Bernoulli random variables formed from the upper tail of the distribution of the index $N$:
 \[Y_i = \left\{
           \begin{array}{rcl}
           1 &~~~& \text{if } N \ge i \\
           0 && \text{otherwise}
           \end{array}
           \right. \]
The Bernoulli sum $X = \sum_{i=1}^\infty Y_i$ is the number of coin tosses ($Z_i$s) that have occurred when the first head ($Z_i =1$) appears in the sequence.  While $X = N$, each provides a different way of looking at the same random variable.

Given the Bernoulli sequence $(Y_{i})_{i\geq1}$ and any $m-$set of indices $\{i_1,\ldots, i_m\}\subseteq \Nset{n}$,  the joint probability of $m$ $Y_i$s  can be written as
\begin{align*}
Pr(Y_{i_1} = 1, \ldots, Y_{i_m} = 1) 
   &= Pr(N \geq i_1, \ldots, N \geq i_m)\\
   &= Pr(N \geq \max{i_1,\ldots, i_m}).
\end{align*}
Then by Proposition \ref{prop:moment_inf_expansion}
\begin{align*}
  E(X^k) &= 
			\sum_{m=1}^{k}S(k,m) \sum_{\{i_1,\ldots, i_m\}\subset \Naturals} 
				E(Y_{i_1} \dotsb Y_{i_m})\\
		 &= \sum_{m=1}^{k}S(k,m) \sum_{\{i_1,\ldots, i_m\}\subset \Naturals} 
		 		Pr(N \geq \max{i_1,\ldots, i_m})\\
		 &= \sum_{m=1}^{k}S(k,m) 
		 		\sum_{M\geq 1} \sum_{\{i_1,\ldots, i_{m-1}\}\subseteq \Nset{M-1}}
		 			Pr(N \ge M)
					~~~\text{(}  M \text{ being the max index)}\\
		 &= \sum_{m=1}^{k}S(k,m) 
		 		\sum_{M\geq 1} \binom{M-1}{m-1} 
				[Pr(N > M) + Pr(N = M)]
		 			\\
		 &=\sum_{m=1}^{k}S(k,m) 
		 		\sum_{M\geq 1} 
		 			\binom{M-1}{m-1} [(1-p)^{M} + (1-p)^{M-1}p]\\
		 &=\sum_{m=1}^{k}S(k,m) 
		 		\sum_{M\geq 1} 
		 			\binom{M-1}{m-1} (1-p)^{M-1}. \\
\end{align*}
For an indeterminate $y$ and $k\in \Naturals$,
$$\sum_{n\geq 0}\binom{n}{k}y^n = \frac{y^k}{(1-y)^{k+1}} $$
from which it follows that
\[\sum_{M\geq 1} 
	\binom{M-1}{m-1} (1-p)^{M-1}
		 =  \frac{(1-p)^{m-1}}{p^m}
		 			\]
giving the expression for the $k$th moment to be		 			
\[
E(X^k) = \sum_{m=1}^{k}S(k,m) 
		 		 \frac{(1-p)^{m-1}}{p^m}.
\]
		 		
Similarly, the $k$th factorial moment has the simpler expression
\[
  E(\ffact{X}{k}) = \frac{(1-p)^{k-1}}{p^{k}}.
\]
We note that since the probability generating function of $X$ has the closed form expression 
\begin{align*}
  G_X(s) = \frac{p}{1-s(1-p)},
\end{align*}
by Equation \ref{eq:fmgf}, the factorial moment generating function is given by
\begin{align*}
  H_X(s) &= \frac{p}{1-(1+s)(1-p)}\\
  &= p\sum_{\ell \geq 0} [(1+s)(1-p)]^{\ell}.
\end{align*}

\subsection{Poisson distribution}
In Section \ref{sec:binom_to_poisson} the tidy expression 
\[E(N^k) = \sum_{m=1}^{k} S_2(k,m)\lambda^m\]
for the $k$th moment of  $N \follows Poisson(\lambda)$ appeared. The proof of this result given by 
 \citeA{riordan1937moment} is recursive. It can also be proved by direct application of Proposition \ref{prop:discrete_moments} as follows:
\begin{align*}
  E(N^k) &= \sum_{m=1}^{k} 
  				S(k,m)
  					 \sum_{M=m}^{\infty}
  					 	\binom{M-1}{m-1}
  					 		Pr(N \geq M)\\
  		 &= \sum_{m=1}^{k}
  		 		S(k,m)
  		 			\sum_{M=m}^{\infty}
  					 	\binom{M-1}{m-1}
  					 		\sum_{\ell = M}^{\infty}e^{-\lambda}\frac{\lambda^\ell}{\ell!}\\
  		 &= e^{-\lambda}\sum_{m=1}^{k} ~
  		 		S(k,m)~
  		 				\sum_{\ell=m}^{\infty}
  					 			\frac{\lambda^\ell}{\ell!} ~
  					 				\sum_{M = 1}^{\ell}\binom{M-1}{m-1}
									~~~ \text{by reorganizing the sums} \\					 		
  		 &= e^{-\lambda}\sum_{m=1}^{k}
  		 		S_2(k,m)m! ~
  		 				\sum_{\ell=m}^{\infty}
  					 			\frac{\lambda^\ell}{\ell!}
  					 				\binom{\ell}{m}\\
  					 				\intertext{as $\sum_{M = 1}^\ell\binom{M-1}{m-1}=\binom{\ell}{m}$,}
		 &= e^{-\lambda}\sum_{m=1}^{k} 
  		 		S_2(k,m)\lambda^m ~
  		 				\sum_{\ell=m}^{\infty}
  					 			\frac{\lambda^{(\ell-m)}}{(\ell-m)!}\\
		 &= e^{-\lambda}\sum_{m=1}^{k}
  		 		S_2(k,m)\lambda^m
  		 				e^{\lambda}\\
		 &= \sum_{m=1}^{k}
  		 		S_2(k,m)\lambda^m. 		 				  					 			   					 							
\end{align*}

Following the same route as for $E(N^k)$, the $k$th factorial moment of $N \follows Poisson(\lambda)$ has the even simpler expression:
\[ E(\ffact{N}{k}) = \lambda^k. \]


\subsection{Ideal soliton distribution}
For an integer $r$ with $r\geq 2$, we say that $N$ follows the ideal soliton distribution, $soliton(r)$, when
\begin{align*}
  Pr(N = 1) &= \frac{1}{r},\\
  Pr(N = i) &= \frac{1}{i(i-1)},
\end{align*}
for $i\in \{2,\ldots, r\}$ and is zero otherwise. It can be shown by induction that 
$$\sum_{i = 2}^{\ell} \frac{1}{i(i-1)} = \frac{\ell-1}{\ell}.$$
From this it follows that
\[Pr(N \leq \ell) 
   = \frac{1}{r} + \frac{\ell - 1}{\ell}  ~~~\text{for } \ell = 1, 2, \ldots r
\]
and 
\[ 
Pr(N \geq \ell)   = \left\{
                                \begin{array}{ccl}
                                   1 & ~~~& \text{when } \ell = 1 \\
                                   && \\
                                   \frac{r-1}{r} - \frac{\ell - 2}{\ell - 1} && \ell = 2, 3, \dots, r.
                                \end{array}
                                \right.
\]
To derive the moments and factorial moments, we need the following lemma.

\begin{lemma}
	\label{lemma:falling_fct_sums}
For any $m, r \in \Naturals$,
$$\sum_{M = 0}^r M\ffact{M}{m} = (m+1)!\binom{r+1}{r-m-1} + m\cdot m!\binom{r+1}{r-m}.$$
\end{lemma}

\begin{proof}
First, an expression for the ordinary generating function $f(x) = \sum_{i\geq 0}i\ffact{i}{m-2}~x^j$ corresponding to the sequence $(i\ffact{i}{m-2})_{i\geq 0}$ is determined.  Multiplying by $\frac{1}{1-x}$ then gives the generating series for $(\sum_{i=0}^{r}i\ffact{i}{m-2})_{r\geq 0}$.

Since $\sum_{i\geq 0}x^i = \frac{1}{1-x}$, differentiating with respect to $x$ and then multiplying by $x$ gives
\begin{align*}
  \sum_{i\geq 0}ix^i &= \frac{x}{(1-x)^2}.
\end{align*}
Differentiating with respect to $x$, $(m-2)$ times, the left hand side becomes
$$ \sum_{i\geq 0}i \times i(i-1)\dotsb (i-m-1)~x^{i-m-2} = \sum_{i\geq 0}i\ffact{i}{m-2} ~x^{i-m-2}.$$
Applying the general Leibniz rule, the same derivative of right hand side is 
$$\frac{x(m-1)!}{(1-x)^m} + (m-2)\frac{(m-2)!}{(1-x)^{m-1}}.$$
Multiplying both sides by $\dfrac{x^{m+2}}{(1-x)}$ gives
\begin{align*}
  \sum_{r\geq 0}\sum_{i=0}^{r}i\ffact{i}{m-2}x^{r} &= \frac{x^{m+3}(m-1)!}{(1-x)^{m+1}} + (m-2)\frac{x^{m+2}(m-2)!}{(1-x)^{m}}.
\end{align*}
Since $$\frac{1}{(1-x)^{k}} = \sum_{n\geq 0} \binom{n+k-1}{n} x^n,$$
the claim immediately follows.
\end{proof}

\begin{proposition}
	\label{prop:soliton_moments}
	Let $N$ follow the soliton($r$) distribution with $r\geq 2$. For $k\geq 1$, the moments of $N$ are given by 
	\begin{align*}
		  E(N^k) &= H_r + \sum_{m=2}^{k}S(k,m)\left[\frac{r-1}{r}\binom{r}{m} - \binom{r-3}{m} - (m-2)^2\binom{r-3}{m-1} \right]
	\end{align*}
	For $k\geq 2$, the $k-$th factorial moment of $N$ is
\[  E(\ffact{N}{k}) = k!\left[\frac{r-1}{r}\binom{r}{k} - \binom{r-3}{k} - (k-2)^2\binom{r-3}{k-1} \right]. \]
\end{proposition}

\begin{proof}
Once more, by Proposition \ref{prop:bern_moments} and Proposition \ref{prop:bern_choose_function}, we need to evaluate sums of the form 
\[ \sum_{M = m}^r
  	\binom{M-1}{m-1}
  		Pr(N \geq M). \]
If $m = 1$, this becomes
\begin{align*}
  \sum_{M = m}^r
  	\binom{M-1}{m-1}
  		Pr(N \geq M) &= \sum_{M=1}^{r} Pr(N \geq 1)\\
  		&= \sum_{M=1}^{r} \sum_{\ell = M}^{r} Pr(N = \ell)\\
  		&= \sum_{\ell=1}^{r}  \sum_{M=1}^{\ell}Pr(N = \ell)\\
  		&= \sum_{\ell=1}^{r} \ell ~ Pr(N = \ell) \\
  		&= \frac{1}{r} + \sum_{j=2}^{r} \frac{j}{j(j-1)} = H_r.
\end{align*}
Otherwise, for $m\geq 2$,   		
\begin{align*}
  \sum_{M = m}^r
  	\binom{M-1}{m-1}
  		Pr(X \geq M) &= 
  \sum_{M = m}^r
  	\binom{M-1}{m-1}
  		\left[ \frac{r-1}{r} - 
  			\frac{M-2}{M-1}\right]\\
  &= \frac{r-1}{r} 
  \sum_{M=m}^{r} 
  	\binom{M-1}{m-1} -
  	\sum_{M=m}^{r}\binom{M-1}{m-1}\frac{M-2}{M-1}.
\end{align*}
By straightforward algebraic manipulation,
$$\sum_{M=m}^r\binom{M-1}{m-1}\frac{M-2}{M-1} = \frac{1}{(m-1)!}\sum_{M=m}^r(M-2)\ffact{M-2}{m-2}.$$
By replacing $m$ and $r$ by $m-2$ and $r-2$, respectively, in Lemma \ref{lemma:falling_fct_sums} and dividing by $(m-1)!$ it follows that
\begin{align*}
  \frac{1}{(m-1)!}\sum_{M=m}^r(M-2)\ffact{M-2}{m-2} &= 
  \frac{1}{(m-1)!}\sum_{M=m}^{r-2}(M)\ffact{M}{m-2}\\ 
  &= \binom{r-1}{r-m-1} + (m-2)^2\binom{r-1}{r-m}.
\end{align*}
The result now follows from a straightforward application of Propositions \ref{prop:bern_choose_function} and \ref{prop:bern_moments}.
\end{proof}

\subsection{Benford distribution}
Benford's distribution encapsulates the notion that in many real world settings, the leading digits in a numerical data set are more likely to be small. In particular, we say that $D$ follows Benford's distribution if for a digit $d\in\{1,\ldots, 9\}$, 
$$Pr(D =d ) = \log_{10}(d+1) - \log_{10}(d).$$
Due to its telescoping nature, the complementary CDF of $D$ is given by
$$Pr(D \geq d) = 1- \log_{10}(d).$$
Therefore, by Proposition \ref{prop:discrete_moments}, the $k-$th moment of $D$ is
\begin{align*}
  E(D^k) &= \sum_{m=1}^{k}S(k,m)
  			\sum_{M = m}^{9}\binom{M-1}{m-1}
  				\left(1-\log_{10}(M) \right)\\
  		 &= \sum_{m=1}^{k}S(k,m)
  		 	\left[ 
  		 		\binom{9}{m} - 
  		 			\sum_{M=m}^{9}
  		 				\binom{M-1}{m-1}
  		 					\log_{10}(M)\right],
\end{align*}
and factorial moments with the form
\begin{align*}
  E(\ffact{D}{k}) &= 
  		 	k!\left[ 
  		 		\binom{9}{k} - 
  		 			\sum_{M=k}^{9}
  		 				\binom{M-1}{k-1}
  		 					\log_{10}(M)\right],
\end{align*}
for $k \leq 9$. Of course, the results above extend to any general base $b$ by noting that $D_b \follows$ Benford$(b)$ satisfies
\begin{align*}
  \sum_{M=m}^{b-1} \binom{M-1}{m-1} Pr(D_b \geq M)&= 
  		 	\left[ 
  		 		\binom{b-1}{m} - 
  		 			\sum_{M=k}^{b-1}
  		 				\binom{M-1}{k-1}
  		 					\log_{b}(M)\right].
\end{align*}

\section{Concluding remarks}
\label{section_bernoulli_summable_conclusion}
A novel multinomial theorem for commutative idempotents (Proposition \ref{prop:multinom_specialization}) led to new general expressions for the moments (including central and factorial) of a Bernoulli sum (e.g., Propositions \ref{prop:bern_moments} to \ref{prop:moment_inf_expansion}) 
as well as corresponding generating functions.  
The general expressions depend on the determination of the expected product of subsets of the Bernoulli random variables. 
By evaluating these in particular cases a number of new expressions for moments and generating functions of many common distributions and classic problems.

The success of the approach in these examples mark it as potentially fruitful in more novel distributions and problems where this expectation might be more readily available.
To that end, the representation of $\binom{X}{m}$ for random count $X$ and fixed $m$ (Proposition \ref{prop:bern_choose_function}), and of $X!$ (Proposition \ref{prop:expected_factorial} ), as the product of Bernoullis may also be more generally useful.

In other instances, the general representation of the various moments for a count variable $N$ expressed in terms of the upper probability of that $N$ (Proposition \ref{prop:discrete_moments}) may be valuable in yet other problems, as shown in the examples of Section \ref{sec:counts_generally}.
 
The framework of Bernoulli sum random variables appears a viable tool for problems involving count data. 

\bibliography{proposal_bibliography.bib}

\begin{thebibliography}{}

\bibitem [\protect \citeauthoryear {%
Chakraborti%
, Jardim%
\BCBL {}\ \BBA {} Epprecht%
}{%
Chakraborti%
\ \protect \BOthers {.}}{%
{\protect \APACyear {2019}}%
}]{%
Chakra2019}
\APACinsertmetastar {%
Chakra2019}%
\begin{APACrefauthors}%
Chakraborti, S.%
, Jardim, F.%
\BCBL {}\ \BBA {} Epprecht, E.%
\end{APACrefauthors}%
\unskip\
\newblock
\APACrefYearMonthDay{2019}{}{}.
\newblock
{\BBOQ}\APACrefatitle {Higher-Order Moments Using the Survival Function: The
  Alternative Expectation Formula} {Higher-order moments using the survival
  function: The alternative expectation formula}.{\BBCQ}
\newblock
\APACjournalVolNumPages{The American Statistician}{73}{2}{191-194}.
\newblock
\begin{APACrefDOI} \doi{10.1080/00031305.2017.1356374} \end{APACrefDOI}
\PrintBackRefs{\CurrentBib}

\bibitem [\protect \citeauthoryear {%
de Montmort%
}{%
de Montmort%
}{%
{\protect \APACyear {1713}}%
}]{%
montmort1713essay}
\APACinsertmetastar {%
montmort1713essay}%
\begin{APACrefauthors}%
de Montmort, P\BPBI R.%
\end{APACrefauthors}%
\unskip\
\newblock
\APACrefYear{1713}.
\newblock
\APACrefbtitle {Essay d'analyse sur les jeux de hazard...} {Essay d'analyse sur
  les jeux de hazard...}
\newblock
\APACaddressPublisher{}{J. Quillau}.
\PrintBackRefs{\CurrentBib}

\bibitem [\protect \citeauthoryear {%
Fr{\'e}chet%
}{%
Fr{\'e}chet%
}{%
{\protect \APACyear {1943}}%
}]{%
frechet1943extension}
\APACinsertmetastar {%
frechet1943extension}%
\begin{APACrefauthors}%
Fr{\'e}chet, M.%
\end{APACrefauthors}%
\unskip\
\newblock
\APACrefYearMonthDay{1943}{}{}.
\newblock
{\BBOQ}\APACrefatitle {Sur l'extension de certaines {\'e}valuations
  statistiques au cas de petits {\'e}chantillons} {Sur l'extension de certaines
  {\'e}valuations statistiques au cas de petits {\'e}chantillons}.{\BBCQ}
\newblock
\APACjournalVolNumPages{Revue de l'Institut International de
  Statistique}{}{}{182--205}.
\PrintBackRefs{\CurrentBib}

\bibitem [\protect \citeauthoryear {%
Goulden%
\ \BBA {} Jackson%
}{%
Goulden%
\ \BBA {} Jackson%
}{%
{\protect \APACyear {2004}}%
}]{%
goulden2004combinatorial}
\APACinsertmetastar {%
goulden2004combinatorial}%
\begin{APACrefauthors}%
Goulden, I\BPBI P.%
\BCBT {}\ \BBA {} Jackson, D\BPBI M.%
\end{APACrefauthors}%
\unskip\
\newblock
\APACrefYear{2004}.
\newblock
\APACrefbtitle {{Combinatorial Enumeration}} {{Combinatorial Enumeration}}.
\newblock
\APACaddressPublisher{}{Dover Publications}.
\PrintBackRefs{\CurrentBib}

\bibitem [\protect \citeauthoryear {%
Iyer%
}{%
Iyer%
}{%
{\protect \APACyear {1958}}%
}]{%
iyer1958theorem}
\APACinsertmetastar {%
iyer1958theorem}%
\begin{APACrefauthors}%
Iyer, P\BPBI K.%
\end{APACrefauthors}%
\unskip\
\newblock
\APACrefYearMonthDay{1958}{}{}.
\newblock
{\BBOQ}\APACrefatitle {A theorem on factorial moments and its applications} {A
  theorem on factorial moments and its applications}.{\BBCQ}
\newblock
\APACjournalVolNumPages{The Annals of Mathematical
  Statistics}{29}{1}{254--261}.
\PrintBackRefs{\CurrentBib}

\bibitem [\protect \citeauthoryear {%
Johnson%
, Kemp%
\BCBL {}\ \BBA {} Kotz%
}{%
Johnson%
\ \protect \BOthers {.}}{%
{\protect \APACyear {2005}}%
}]{%
johnson2005univariate}
\APACinsertmetastar {%
johnson2005univariate}%
\begin{APACrefauthors}%
Johnson, N\BPBI L.%
, Kemp, A\BPBI W.%
\BCBL {}\ \BBA {} Kotz, S.%
\end{APACrefauthors}%
\unskip\
\newblock
\APACrefYear{2005}.
\newblock
\APACrefbtitle {{Univariate Discrete Distributions}} {{Univariate Discrete
  Distributions}}\ (\BVOL~444).
\newblock
\APACaddressPublisher{}{John Wiley \& Sons}.
\PrintBackRefs{\CurrentBib}

\bibitem [\protect \citeauthoryear {%
Knoblauch%
}{%
Knoblauch%
}{%
{\protect \APACyear {2008}}%
}]{%
knoblauch2008closed}
\APACinsertmetastar {%
knoblauch2008closed}%
\begin{APACrefauthors}%
Knoblauch, A.%
\end{APACrefauthors}%
\unskip\
\newblock
\APACrefYearMonthDay{2008}{}{}.
\newblock
{\BBOQ}\APACrefatitle {Closed-form expressions for the moments of the binomial
  probability distribution} {Closed-form expressions for the moments of the
  binomial probability distribution}.{\BBCQ}
\newblock
\APACjournalVolNumPages{SIAM Journal on Applied Mathematics}{69}{1}{197--204}.
\PrintBackRefs{\CurrentBib}

\bibitem [\protect \citeauthoryear {%
Riordan%
}{%
Riordan%
}{%
{\protect \APACyear {1937}}%
}]{%
riordan1937moment}
\APACinsertmetastar {%
riordan1937moment}%
\begin{APACrefauthors}%
Riordan, J.%
\end{APACrefauthors}%
\unskip\
\newblock
\APACrefYearMonthDay{1937}{}{}.
\newblock
{\BBOQ}\APACrefatitle {Moment recurrence relations for binomial, Poisson and
  hypergeometric frequency distributions} {Moment recurrence relations for
  binomial, poisson and hypergeometric frequency distributions}.{\BBCQ}
\newblock
\APACjournalVolNumPages{The Annals of Mathematical Statistics}{8}{2}{103--111}.
\PrintBackRefs{\CurrentBib}

\bibitem [\protect \citeauthoryear {%
Shah%
}{%
Shah%
}{%
{\protect \APACyear {1973}}%
}]{%
shah1973distribution}
\APACinsertmetastar {%
shah1973distribution}%
\begin{APACrefauthors}%
Shah, B.%
\end{APACrefauthors}%
\unskip\
\newblock
\APACrefYearMonthDay{1973}{}{}.
\newblock
{\BBOQ}\APACrefatitle {Distribution of sum of independent integer valued
  random-variables} {Distribution of sum of independent integer valued
  random-variables}.{\BBCQ}
\newblock
\APACjournalVolNumPages{The American Statistician}{27}{3}{123--124}.
\PrintBackRefs{\CurrentBib}

\bibitem [\protect \citeauthoryear {%
Shmueli%
, Minka%
, Kadane%
, Borle%
\BCBL {}\ \BBA {} Boatwright%
}{%
Shmueli%
\ \protect \BOthers {.}}{%
{\protect \APACyear {2005}}%
}]{%
shmueli2005useful}
\APACinsertmetastar {%
shmueli2005useful}%
\begin{APACrefauthors}%
Shmueli, G.%
, Minka, T\BPBI P.%
, Kadane, J\BPBI B.%
, Borle, S.%
\BCBL {}\ \BBA {} Boatwright, P.%
\end{APACrefauthors}%
\unskip\
\newblock
\APACrefYearMonthDay{2005}{}{}.
\newblock
{\BBOQ}\APACrefatitle {{A useful distribution for fitting discrete data:
  revival of the Conway--Maxwell--Poisson distribution}} {{A useful
  distribution for fitting discrete data: revival of the
  Conway--Maxwell--Poisson distribution}}.{\BBCQ}
\newblock
\APACjournalVolNumPages{Journal of the Royal Statistical Society: Series C
  (Applied Statistics)}{54}{1}{127--142}.
\PrintBackRefs{\CurrentBib}

\bibitem [\protect \citeauthoryear {%
Stanley%
}{%
Stanley%
}{%
{\protect \APACyear {2011}}%
}]{%
stanley2011enumerative}
\APACinsertmetastar {%
stanley2011enumerative}%
\begin{APACrefauthors}%
Stanley, R\BPBI P.%
\end{APACrefauthors}%
\unskip\
\newblock
\APACrefYear{2011}.
\newblock
\APACrefbtitle {{Enumerative Combinatorics, Volume I}} {{Enumerative
  Combinatorics, Volume I}}\ (\PrintOrdinal{2nd}\ \BEd).
\PrintBackRefs{\CurrentBib}

\bibitem [\protect \citeauthoryear {%
Wilf%
}{%
Wilf%
}{%
{\protect \APACyear {2005}}%
}]{%
wilf2005generatingfunctionology}
\APACinsertmetastar {%
wilf2005generatingfunctionology}%
\begin{APACrefauthors}%
Wilf, H\BPBI S.%
\end{APACrefauthors}%
\unskip\
\newblock
\APACrefYear{2005}.
\newblock
\APACrefbtitle {generatingfunctionology} {generatingfunctionology}.
\newblock
\APACaddressPublisher{}{A K Peters/CRC Press}.
\PrintBackRefs{\CurrentBib}

\end{thebibliography}

\end{document}